\documentclass[12pt,leqno]{article}

\newcommand\qed{\hfill$\sqcap\kern-8.0pt\hbox{$\sqcup$}$}
\newcommand\NN{\hbox{I\kern-.2em\hbox{N}}}
\newcommand\RR{\hbox{I\kern-.2em\hbox{R}}}
\newcommand\sRR{{\sl \hbox{I\kern-.2em\hbox{R}}}}
\newcommand{\PP}{{\bf P}^k}

\usepackage{amsmath}
\input{amssym}
\usepackage{amscd}
\usepackage[latin1]{inputenc}
\usepackage[french]{babel}
\newcommand\QQ{\hbox{I\kern-.53em\hbox{Q}}}
\newcommand\ZZ{{{\rm Z}\kern-.28em{\rm Z}}}
\newtheorem{theo}{Th\'eor\`eme}
\newtheorem{prop}{Proposition}
\newtheorem{lem}{Lemme}
\newtheorem{rem}{Remarque}

\newtheorem{defi}{D\'efinition}

\begin{document}

\date{}

\title {Linéarisation d'endomorphismes holomorphes de ${\bf P}^k$
et caractérisation des exemples de Lattès par leur mesure
de Green}
\vskip0.5cm
\author { F. Berteloot et C. Dupont}
\vskip0.5cm
\maketitle

\section{Introduction et résultats}
Les propriétés dynamiques d'un endomorphisme $f$, holomorphe et de degré algébri-que $d\ge 2$ sur l'espace projectif
complexe ${\bf P}^k$, se reflètent sur son courant et sa mesure de Green. Ce sont respectivement 
un $(1,1)$-courant positif fermé $T$ obtenu comme limite de $\frac{1}{d^n}f^{n*} \omega$ où $\omega$
désigne la forme de Fubiny-Study et une mesure de probabilité invariante $\mu$ obtenue comme $k$-ième
puissance extérieure de $T$. Ces objets, introduits par Hubbard-Papadopol \cite{HP} et Fornaess-Sibony
\cite{FS} possèdent de remarquables propriétés ergodiques. Fornaess et Sibony ont montré
que la mesure de Green est mélangeante \cite{FS2}, elle est aussi l'unique mesure d'entropie maximale et ses exposants de
Liapounov sont supérieurs à $\frac{1}{2}Log\;d$ comme l'ont montré Briend et Duval \cite{BD}, \cite{BD2}.\\

La dimension de Hausdorff de $\mu$ ($Dim_H\mu$), définie comme la borne inférieure des dimensions de
Hausdorff des ensembles de $\mu$-mesure pleine, est une caractéristique géométrique importante du système
dynamique $({\bf P}^k, f, d, \mu)$. L'une des premières questions concernant l'estimation de cette
dimension est de déterminer les
systèmes pour lesquels elle est maximale.\\
En dimension $k=1$, Ledrappier a montré que la maximalité de $Dim_H\mu$ 
équivaut à l'absolue
continuité de $\mu$ par rapport à la mesure de Lebesgue \cite{L2}. Un argument de renormalisation 
montre alors que l'endomorphisme
$f$ est un exemple de Lattès, il s'agit d'un cas particulier d'un résultat de Mayer \cite{M}. Notons
aussi le résultat très précis de Zdunik \cite {Z} qui stipule que la dimension
de Hausdorff de $\mu$ coïncide avec celle de son support (l'ensemble de Julia de $f$) si et seulement
si $f$ est un exemple de Lattès, un polyn\^ ome de Tchebychev ou une puissance $z^{\pm d}$.\\                   
En dimension supérieure, il est possible d'adapter le travail de Ledrappier selon
lequel,
pour toute mesure invariante par $f$, l'égalité dans l'inégalité de Margulis-Ruelle force l'absolue
continuité. Ceci fait l'objet de \cite{D} et concerne en particulier les mesures de Green 
d'exposants de Liapounov minimaux. Par ailleurs, la maximalité de $Dim_H\mu$ entra\^ \i ne la minimalité des
exposants,   
comme l'ont implicitement établi Binder et DeMarco \cite{BdM} (voir aussi l'appendice).
Ainsi, une caractérisation précise des endomorphismes dont la mesure est
absolument continue montrerait que, pour un système $({\bf P}^k, f, d, \mu)$  générique,
l'un au moins des exposants de Liapounov de $\mu$ est strictement supérieur à $\frac{1}{2}Log\;d$ 
et $Dim_H\mu$ est strictement
inférieure à $2k$. En outre ceci répondrait à une question posée par Fornaess et Sibony dans \cite{FS4}.
Le principal résultat de cet article fournit une telle caractérisation :  

\begin{theo}\label{TH} 
Pour tout système $({\bf P}^k, f, d, \mu)$, la mesure de Green $\mu$ est non singulière par rapport à
la mesure de Lebesgue si et seulement si $f$ est un endomorphisme de Lattès :
il existe alors un diagramme commutatif 
\begin{equation*}
\begin{CD}
{\bf C}^k @>D>>{\bf C}^k\\
{\sigma}@VVV                          @VV{\sigma}V\\
{\bf P}^k @>f>>                {\bf P}^k
\end{CD}
\end{equation*}    
où $D$ est une application affine de partie
lin\'eaire
$\sqrt d\;U$ ($U$ unitaire)et $\sigma$ un rev\^etement ramifi\'e
sur les fibres duquel un groupe cristallographique complexe  
agit transitivement.
\end {theo}

Il est crucial, dans notre contexte, que la mesure de Green provienne d'un courant ($\mu=T^k$) car
l'invariance du courant de Green ($f^*T=dT$) recèle plus d'informations géométriques que celle
de la mesure. En particulier, il a été démontré dans \cite{BL} que $f$ est un endomorphisme de Lattès dès que
son courant de Green est lisse et strictement positif sur un ouvert de ${\bf P}^k$. La preuve du
théorème \ref{TH} revient donc à déduire la régularité de $T$ de celle de $\mu$. Nous utiliserons pour cela
des arguments de renormalisation qui nécessitent de linéariser la suite $(f^n)_n$ le long d'orbites
typiques.\\

La mise au point d'un tel procédé de linéarisation occupe la section 3. Il s'agit, pour des choix 
$\mu$-génériques de $x$, de rendre la suite $(f^n)_n$ normale en $x$ en la 
précomposant par
des contractions équivalentes aux applications linéaires
tangentes
inverses $\big(\tau_x \circ (d_0f^{n}_x)^{-1}\big)$. A cet effet, nous estimons
précisément les erreurs cumulées lorsque l'on remplace $f$ par sa différentielle le long d'une orbite
(cf. Proposition \ref{pr3.3}). Outre la stricte positivité des exposants de Liapounov 
$\lambda_1 \le ...\le \lambda_k$ du système $({\bf P}^k, f, d, \mu)$, ceci requiert l'hypothèse        
$\lambda_k<2\lambda_1$. Une fois acquise la possibilité de linéariser, nous majorons la norme des
différentielles $(d_0f^{n}_x)^{-1}$ en reprenant, dans ce contexte,
la méthode pluripotentialiste de Briend et Duval. Nous obtenons finalement une minoration de la masse des
points $x$ où $(f^n)_n$ est linéarisable et les normes $\vert\vert(d_0f^{n}_x)^{-1}\vert\vert$
convenablement majorées :

\begin{theo}\label{THB} Soit $({\bf P}^k, f, d, \mu)$ un système tel que  
$\lambda_k < 2 \lambda_1$. Pour $\tau >0$, $\rho\in]0,1]$ et $n\in{\bf N}$, soit 
${\cal L}{\cal B}_n(\rho,\tau)$ l'ensemble des points $x\in{\bf P}^k$ tels que 
$f^n_x \circ (d_0f^n_x)^{-1}$ soit injective de $B(0,\rho)$ dans $B(0,R_0)$
et $\vert\vert (d_0f^n_x)^{-1} \vert\vert
\le \tau d^{-\frac{n}{2}}$.\\
Alors il existe $\alpha: ]0,1] \to {\bf
R}^+$ et $C>0$ ne dépendants que de $f$
tels que $\lim_{\rho \to 0} \alpha(\rho) = 1$ et 
$\liminf_n \mu[{\cal L}{\cal B}_n(\rho,\tau)]\ge\alpha(\rho)-\frac{C}{\tau^2\rho^2}$. 

\end{theo}

Il est facile de d'établir une version plus maniable du théorème \ref{THB}. En outre, gr\^ ace à
la propriété de mélange, on peut assujettir une sous-suite de l'orbite $(f^n(x))_n$ 
à ne pas s'échapper d'un borélien
prescrit :

\begin{theo}\label{THL} Si les exposants de Liapounov de $({\bf P}^k, f, d, \mu)$ sont tels que 
$\lambda_k < 2 \lambda_1$, alors pour tout borélien $B$ chargé par $\mu$, on peut 
trouver $\tilde B \subset
B$ de masse arbitrairement proche de $\mu(B)^2$ et $\tau_0>0$ tels que pour tout point
$x \in \tilde B$ il existe une suite extraite $(f^{n_j})_j$ ainsi qu'un réel $\nu(x)>0$ vérifiant :
\begin{itemize}
\item[1)] $f^{n_j}(x) \in B$ pour tout $j\in {\bf N}$
\item[2)] $f^{n_j} \circ \tau_x \circ (d_0f^{n_j}_x)^{-1}$ converge uniformément vers un biholomorphisme
sur $B(0,\nu(x))$
\item[3)] $\vert\vert (d_0f^{n_j}_x)^{-1}\vert\vert \le \tau_0 d^{-\frac{n_j}{2}}$.

\end{itemize}
\end{theo}

La section 4 concerne
la 
preuve du
théorème \ref{TH} proprement dite. La régularité de $T$ 
se déduit, par des arguments de
renormalisation, des relations $f^{n*}T=d^nT$ pourvu que la suite $(f^n)_n$ soit assimilable à une suite 
d'homothéties $(d^{-\frac{n}{2}}Id)_n$ (cf. Lemme \ref{lem4.1}). Il s'agit donc de s'assurer que, lorsque
$\mu$ est absolument continue, les différentielles $(d_0f^{n_j}_x)^{-1}$ intervenant dans le théorème
\ref{THL} sont équivalentes à des homothéties de rapport $d^{-\frac{n_j}{2}}$. En d'autres termes,
il faut contr\^ oler les distorsions des ellipsoïdes $(d_0f^{n_j}_x)^{-1}\big[B(0,1)\big]$. Or, d'après la
dernière assertion du théorème \ref{THL} la
taille de ces ellipsoïdes est au plus de l'ordre de $d^{-\frac{n_j}{2}}$ et il convient donc d'en minorer
convenablement le volume. Ceci résulte de l'absolue continuité de $\mu$ et des relations $f^{n*}\mu=d^{kn}
\mu$
(cf. Proposition \ref{pr4.1}).\\

Signalons pour finir que les exemples de Lattès interviennent naturellement dans d'autres problèmes
(voir par exemple \cite{Du1} et \cite{DS}).  

\section{Préliminaires.}

Dans cette section, nous résumons les principaux outils et résultats utilisés dans la suite et fixons
quelques notations.\\

$\bullet$ L'espace projectif complexe ${\bf P}^k$ est muni d'une structure de variété hermitienne induite par
la forme de Fubini-Study $\omega$. On construit une famille de cartes holomorphes 
$\big(\tau_x\big)_{x\in{\bf P}^k}$ telle que :

\begin{itemize} 

\item[1.] $\tau_x :{\bf C}^k \to {\bf P}^k$ est un biholomorphisme sur son image et $\tau_x(0)=x$ 
\item[2.] $(\tau_x^*\omega)_0=\frac{i}{2}\sum_{j=1,k} dz_j\wedge d\bar z_j$\\
\end{itemize}

Cette famille est obtenue en explicitant une telle carte pour un point base $x_0\in {\bf P}^k$ puis en la
propageant à ${\bf P}^k$ par l'action transitive de ${\bf U}_{k+1}({\bf C})$. Ce faisant, on obtient 
plut\^ ot une classe de cartes en $x$ car $\tau_x$ est
définie à un élément du sous-groupe d'isotropie de $x_0$ près.
Cette ambiguïté pourra cependant \^ etre ignorée puisque ${\bf U}_{k+1}({\bf C})$ est compact; les
affirmations faisant intervenir $\tau_x$ devront \^ etre comprises comme valables pour tous les éléments 
de la
classe de cartes en $x$.\\
On peut aussi, localement, faire un choix différentiable de $\tau_x$ et en particulier s'assurer de la propriété
suivante :
    
\begin{itemize}
\item[3.]  
$\tau_{x_0}^{-1}\circ\tau_x - \tau_{x_0}^{-1}(x)$ converge vers l'identité en topologie ${\cal C}^{\infty}$
lorsque $x \to x_0$.
\end{itemize}

$\bullet$ L'extension naturelle $(\widehat{{\bf P}^k}, \hat f, \hat\mu)$ est un système dynamique inversible
associé au système $({\bf P}^k, f, d, \mu)$ de la fa\c con suivante :\\
$\widehat{{\bf P}^k}:=\{\hat x:=(x_n)_{n\in {\bf Z}} / f(x_n)=x_{n+1}\}$ est muni de la topologie et de la
tribu produit.\\
On note $\pi_0 : \widehat{{\bf P}^k}\to {{\bf P}^k}$ la projection définie par $\pi_0(\hat x)=x_0$ de sorte
que $\pi_0 \circ f = \hat f \circ \pi_0$ où $\hat f$ désigne le décalage à droite sur $\widehat{{\bf P}^k}$.
La mesure $\hat {\mu}$ est l'unique mesure de probabilité invariante par $\hat f$ sur $\widehat{{\bf P}^k}$
telle que $\pi_{0*}\hat\mu=\mu$. Elle hérite de $\mu$ le caractère mélangeant.\\

On notera ${\hat f}^{-n}$ le décalage à gauche itéré $n$ fois.\\

Soit ${\cal C}_f$ l'ensemble des points critiques de $f$, 
on considère alors $$\widehat{X}:=\{\hat x \in \widehat{{\bf P}^k} /x_n \notin {\cal C}_f, \forall n \in
{\bf Z}  \}$$
\noindent cet ensemble est de $\hat\mu$-mesure pleine car $\mu$ ne charge pas les ensembles pluripolaires
(voir \cite{S} Proposition A.6.3).\\

$\bullet$ Nous noterons $B(0,R)$ (resp. $P(0,R)$) la boule euclidienne centrée en $0$ et de rayon $R$
(resp. le polydisque centré en $0$ et de polyrayon $R$) de ${\bf C}^k$. Nous noterons $B(x,s)$ l'image de 
$B(0,s)$ par $\tau_x$\\

$\bullet$ A tout endomorphisme holomorphe $f$ sur ${\bf P}^k$ et tout $x\in{\bf P}^k$ on associe 
les applications suivantes. Elles sont définies sur un voisinage de l'origine de ${\bf C}^k$
dont la taille dépend de $x$ et $n$ :

$$f_x:=\tau_{f(x)}^{-1} \circ f \circ \tau_x$$
$$f_x^n=\tau_{f^n(x)}^{-1} \circ f^n \circ \tau_x=f_{f^n(x)}\circ ...\circ f_x$$

Pour tout $\hat x\in
\widehat{X}$ on définit une application $f^{-n}_{\hat x}$ 
par :

$$f^{-n}_{\hat x}:=f^{-1}_{x_{-n}}\circ ... \circ f^{-1}_{x_{-1}}$$

$f^{-n}_{\hat x}$ est définie sur un voisinage de l'origine dont la taille dépend mesurablement
de $\hat x$ (voir le lemme ci-dessous).\\

$\bullet$ Nous notons $T=T_a+T_s$ la décomposition de Lebesgue d'un $(1,1)$-courant positif $T$.
Cette décomposition est unique et les courants $T_a$, $T_s$ sont positifs. Par contre la fermeture
de $T$ n'implique pas celle de $T_a$ (ou $T_s$).\\   
Nous noterons $\sigma_T$ la mesure trace du courant $T$. On observe que la décomposition de Lebesgue de
$\sigma_T$ est donnée par $\sigma_T=\sigma_{T_a}+\sigma_{T_s}$.\\

$\bullet$ Nous utiliserons le résultat de Briend et Duval sur la minoration optimale 
des exposants de Liapounov :\\

\noindent{\bf Théorème} (Briend-Duval) 
{\it Les exposants de Liapounov du système $({\bf P}^k, f, d, \mu)$
sont supérieurs ou égaux à $\frac{1}{2}Log\;d$}.\\        
    
\noindent ainsi que des propriétés des branches inverses de $f$, mises en évidence dans leur démonstration et qui,
compte tenu de la stricte positivité des exposants peuvent s'exprimer ainsi 
(cf. \cite{BD} ou \cite{Du}) :

\begin{lem}\label{lempr}
Soient $\lambda_1\le ...\le \lambda_k$ les exposants de Liapounov du système $({\bf P}^k, f, d, \mu)$.
Soient $0<\epsilon\le\epsilon_0\ll 1$ et $0<r_0\le R_0\ll 1$. Il existe des fonctions $\rho$, $r$ continues
sur ${\bf P}^k$ et strictement positives hors de ${\cal C}_f$ ainsi que des fonctions mesurables
$\eta : \widehat{Y} \to ]0,r_0]$, $C : \widehat{Y} \to [1,+\infty[$ définies sur un ensemble de
$\hat{\mu}$-mesure pleine $\widehat{Y}$ telles que :\
\begin{itemize} 

\item[1.] $f_x$ est injective sur $B(0,\rho(x))$ et $B(0,r(x)) \subset f_x\big[B(0,\rho(x))\big]$, $\forall
x\in {\bf P}^k \setminus {\cal C}_f$ 
\item[2.] $\lim_n \frac{1}{n}Log\;\rho(x_n) = 0$, $\forall \hat x \in \widehat{Y}$
\item[3.] $f^{-n}_{\hat x}$ est injective sur $B(0,\eta(\hat x))$ pour tout $\hat x\in\widehat{Y}$ et tout
$n$, de plus :  
$$d_0 f^{-n}_{\hat x}\big[B(0,\gamma\eta (\hat x))\big]\subset B\big(0,\gamma
r(x_{-(n+1)})e^{-n(\lambda_1-\epsilon)}\big)$$
pour tout $0<\gamma<1$
\item[4.] $Lip\;f^{-n}_{\hat x} \le C(\hat x) e^{-n(\lambda_1-\frac{\epsilon}{2})}$ sur 
$B(0,\eta (\hat x))$. 

\end{itemize}
     
\end{lem}

\section{Un proc\'ed\'e de linéarisation.}

Dans toute cette partie nous consid\'erons le système dynamique $({\bf P}^k, f, d, \mu)$ et 
adoptons la définition suivante.

\begin{defi}
La suite des itérées 
$(f^n)_n$ est {\it linéarisable en $x \in {\bf P}^k$} si et seulement si il existe $\nu(x)>0$ tel
qu'après une éventuelle extraction, la suite $[f^n \circ \tau_x \circ (d_0f^n_x)^{-1}]_n$ converge
uniformément vers une limite injective sur $B(0,\nu(x))$.
\end{defi}

Nous commencerons par montrer que cette propriété est $\mu$-générique dès que
le spectre de Liapounov du système $({\bf P}^k, f, d, \mu)$ est assez étroit, cela correspond à la
seconde assertion du théorème \ref{THL} lorsque le borélien $B$ est pris égal à ${\bf P}^k$ : 

\begin{theo}\label{th3.1} Si les exposants de Liapounov de $({\bf P}^k, f, d, \mu)$ sont tels que 
$\lambda_k < 2 \lambda_1$ alors $(f^n)_n$
est linéarisable en $\mu$-presque tout point.
\end{theo}    

Fixons $R_0>0$ puis, pour tout $\rho\in ]0,1]$ et tout $n\in {\bf N}$ définissons ${\cal B}_n(\rho)$ par :
$${\cal B}_n(\rho) :=\{x\in{\bf P}^k / f^n_x \circ (d_0f^n_x)^{-1} injective\;de\; B(0,\rho)\; dans \;
B(0,R_0)\}$$

Dans ce qui suit, on pourra diminuer $R_0$ sans affecter la validité des énoncés. La linéarisabilité en
$x$ résulte immédiatement, via le théorème de Montel, de l'appartenance de $x$ à 
$\cup_{0<\rho\le 1}{\cal B}(\rho)$ où 
${\cal B}(\rho):=\limsup_n{\cal B}_n(\rho)$. Ainsi, comme $\mu[{\cal B}(\rho)] \ge 
\limsup_n \mu[{\cal B}_n(\rho)]$, l'énoncé suivant est une version quantifiée du théorème \ref{th3.1} : 

\begin{prop}\label{prop3.2} Si les exposants de Liapounov de $({\bf P}^k, f, d, \mu)$ sont tels que 
$\lambda_k < 2 \lambda_1$, alors il existe $\alpha: ]0,1] \to {\bf
R}^+$ telle que $\lim_{\rho \to 0} \alpha(\rho) = 1$ et $\mu[{\cal B}_n(\rho)]\ge\alpha(\rho)$ pour tout 
$n\in {\bf N}$.
\end{prop}

Pour établir la proposition \ref{prop3.2}, nous nous pla\c cons dans l'extension naturelle $\widehat{{\bf P}^k}$ 
et comparons les branches inverses de $f$ à leurs différentielles le long d'orbites négatives génériques 
${\hat x}_{-}:= (x_{-j})_{j\ge 0}$. Plus précisément, l'invariance de $\hat \mu$ nous permettra de déduire la
proposition \ref{prop3.2} de la proposition technique suivante : 

\begin{prop}\label{pr3.3}
Sous les hypothèses du théorème \ref{th3.1}, pour tout $r_0 \in ]0,R_0]$ il existe des fonctions mesurables 
$\eta, S
: \widehat{Z} \to ]0,r_0]$ définies sur un ensemble $\widehat{Z}$ de $\hat\mu$-mesure pleine telles
que $S\le \eta$ et
$d_0f^{-n}_{\hat x}\big[B(0,S(\hat x))\big] \subset
f^{-n}_{\hat x}\big[B(0,\eta(\hat x))\big]$ pour tout $\hat x \in \widehat{Z}$ et tout $n\in {\bf N}.$   
 
\end{prop}

Voyons comment la proposition \ref{prop3.2} se déduit de la proposition \ref{pr3.3}. Posons $\widehat{\cal S}(\rho):=\{
\hat x \in \widehat{Z} / S(\hat x) \ge \rho   \}$ où $S$ désigne la fonction fournie par la proposition 
\ref{pr3.3} et vérifions les inclusions suivantes :
$$\pi_0 \big[ \hat{f}^{-n}\big(\widehat{\cal S}(\rho) \big) \big] \subset {\cal B}_n(\rho)$$
L'appartenance de ${\hat x}_n:= \hat{f}^{n}(\hat x)$ à $\widehat{\cal S}(\rho)$ signifie que :

$$ d_0f^{-n}_{{\hat x}_n}\big[B(0,\rho)\big] \subset d_0f^{-n}_{{\hat x}_n}\big[B(0,S({\hat x}_n))\big]
\subset
f^{-n}_{{\hat x}_n}\big[B(0,\eta({\hat x}_n)\big]$$
Comme $f^{-n}_{{\hat x}_n}$ est injective sur $B(0,\eta({\hat x}_n))$ d'inverse $f^n_{x_0}$, l'appartenance
de $x$ à ${\cal B}_n(\rho)$ s'obtient en composant les inclusions précédentes par $f^n_{x_0}$ :
$$f^n_{x_0} \circ \big(d_0 f^n_{x_0} \big)^{-1}\big[B(0,\rho)\big] \subset B(0,\eta({\hat x}_n)) \subset
B(0,R_0)$$

Compte tenu de l'invariance de $\hat \mu$ on a $\mu\big[{\cal B}_n(\rho) \big]\ge \hat\mu
\big[ \hat{f}^{-n}\big(\widehat{\cal S}(\rho) \big)\big] = \hat\mu\big[\widehat{\cal
S}(\rho)\big]=:\alpha(\rho)$ et, comme $S$ est $\hat\mu$-presque partout strictement positive, 
$\lim_{\rho\to 0}\alpha(\rho)=1$.    \hfill$\square$\\

La preuve de la proposition \ref{pr3.3} consiste à compenser les erreurs dues à la substitution de        
$d_0 f^{-1}_{x_j}$ à $ f^{-1}_{x_j}$ le long de $\hat x$ en diminuant le rayon $\eta(\hat x)$. Pour que
les compensations, cumulées, fournissent un rayon $S(\hat x) >0$, il faut que les erreurs soient
négligeables devant la plus petite dimension caractéristique de l'ellipsoïde $d_0 f^{-j}_{\hat
 x}\big[B(0,1)\big]$. Il résulte du lemme suivant que tel est le cas lorsque $\lambda_k < 2\lambda_1$.

\begin{lem}\label{lem3.4}
Pour tout système $({\bf P}^k, f, d, \mu)$ et tout $0<\epsilon\ll 1$, on peut trouver $\widehat Z\subset
 \widehat{{\bf P}^k}$ de $\hat\mu$-mesure pleine et des fonctions mesurables $\eta, E, F : 
 \widehat{Z} \to {\bf R}^{+*}$ telles que $0<\eta\le r_0 \le R_0$ et :\\ 
1) $\forall \hat x \in \widehat{Z}, \forall n\in {\bf N}, \forall \gamma\in ]0,1], \forall u \in 
d_0f^{-n}_{\hat x}\big[B(0,\gamma\eta(\hat x))\big]$ : 
$$ \vert\vert \big(d_0f^{-1}_{ x_{-(n+1)}} - f^{-1}_{ x_{-(n+1)}}\big) (u) \vert\vert \le \gamma E(\hat x) 
e^{-2n(\lambda_1-\epsilon)}$$
2) $\forall \hat x \in \widehat{Z}, \forall n\in {\bf N}$ :  
$\vert\vert d_0f^{n+1}_{ x_{-(n+1)}}\vert\vert
\le F(\hat x)e^{n(\lambda_k+\epsilon)}$          
\end{lem}

\underline{Preuve de la Proposition \ref{pr3.3}} : \\
Reprenons les notations du lemme \ref{lem3.4} et définissons sur $\widehat Z$ les fonctions 
mesurables suivantes :

$$\xi_n(\hat x):= Sup\{t\le \eta(\hat x) / d_0f^{-n}_{\hat x}\big[B(0,t)\big]
\subset
f^{-n}_{\hat x}\big[B(0,\eta({\hat x}))\big]\}$$ 

$$n_0(\hat x) := Min \{p \ge 1 / \forall n \ge p : \frac{EF}{\eta}(\hat x)\le e^{n\epsilon} \}$$

$$s(\hat x) := Min\{\xi_n(\hat x) : 0\le n \le n_0(\hat x)\}$$

Posons $\kappa_j:=1-e^{-j(2\lambda_1-\lambda_k-6\epsilon)}$ où $\epsilon$ est assez petit pour que 
$\prod_{j=1}^{\infty}\kappa_j =: \kappa>0$ et définissons $s_n(\hat x)$ par:

$$s_n(\hat x) := s(\hat x)\;\;\; si\; n\le n_0(\hat x)$$
$$ s_n(\hat x) := s(\hat x)\prod_{j=n_0(\hat x)}^{n-1}\kappa_j\;\;\; si\; n\ge n_0(\hat x)+1$$

Pour montrer que la fonction $S(\hat x):=\kappa s(\hat x)$ convient, il suffit alors d'établir les 
inclusions :

$$ (I_n)_{n\ge 0}\;\;\;:\;\;\;d_0f^{-n}_{\hat x}\big[B(0,s_n(\hat x))\big]
\subset
f^{-n}_{\hat x}\big[B(0,\eta({\hat x}))\big]$$   

Par définition de $s_n(\hat x)$, les inclusions $(I_n)$ sont satisfaites lorsque $n\le n_0(\hat x)$.
Supposons
$(I_n)$ vraie pour $n\ge n_0(\hat x)$ et posons $\nu_n:= \big(\frac{E s_n}{\eta}\big)(\hat x)
e^{-2n(\lambda_1-2\epsilon)}$. On a alors : 

\begin{equation}\label{E1}
 s_{n+1}\le s_n - \vert\vert\big( d_0f^{-(n+1)}_{ \hat x}\big)^{-1}\vert\vert \nu_n
\end{equation}
En effet :
$$s_{n+1}=s_n\kappa_n=s_n\big(1-e^{-n(2\lambda_1-\lambda_k -6\epsilon)}\big)\le
s_n\big(1-\frac{EF}{\eta}(\hat x)e^{-n(2\lambda_1-\lambda_k -5\epsilon)}\big)$$
$$\le
s_n - \vert\vert d_0f^{n+1}_{ x_{-(n+1)}}\vert\vert\nu_n = 
s_n - \vert\vert\big( d_0f^{-(n+1)}_{ \hat x}\big)^{-1}\vert\vert \nu_n$$

la première majoration résultant de la définition de $n_0(\hat x)$ et la seconde de l'assertion (2) du
lemme \ref{lem3.4}.\\
Désignons par $\Lambda$ la frontière de $d_0f^{-(n+1)}_{\hat x}\big[B(0,s_n)\big]$, on verifie alors aisément
que l'inégalité (\ref{E1})
se traduit par 
\begin{equation}\label{E2}
d_0f^{-(n+1)}_{\hat x}\big[B(0,s_{n+1})\big] \subset
d_0f^{-(n+1)}_{\hat x}\big[B(0,s_n)\big] \backslash \bigcup_{p\in \Lambda}B(p,\nu_n)
\end{equation}
Par ailleurs, la première assertion du lemme \ref{lem3.4} (où l'on prend $\gamma=\frac{s_n}{\eta}$)
stipule que $f^{-1}_{ x_{-(n+1)}}$ diffère d'au plus $\nu_n$ de sa différentielle sur 
$d_0f^{-n}_{\hat x}\big[B(0,s_{n})\big]$, il s'ensuit que 
\begin{equation}\label{E3}
d_0f^{-(n+1)}_{\hat x}\big[B(0,s_n)\big] \backslash \bigcup_{p\in \Lambda}B(p,\nu_n)
\subset f^{-1}_{ x_{-(n+1)}} \circ d_0f^{-n}_{\hat x}\big[B(0,s_{n})\big]  
\end{equation}
 
Observons finalement que l'inclusion $(I_n)$, composée par $f^{-1}_{ x_{-(n+1)}}$, s'écrit

\begin{equation}\label{E4}
f^{-1}_{ x_{-(n+1)}} \circ d_0f^{-n}_{\hat x}\big[B(0,s_{n})\big]
\subset
f^{-(n+1)}_{\hat x}\big[B(0,\eta({\hat x})\big]
\end{equation}
Les inclusions (\ref{E2}), (\ref{E3}) et (\ref{E4}) encha\^\i nées donnent $(I_{n+1})$
\hfill$\square$    
\\

\underline{Preuve du lemme \ref{lem3.4}} :\\

D'après la première assertion du lemme \ref{lempr}, l'application $f_{ x_{-(n+1)}}$ est inversible sur 
$B(0,r)$ où $r:=r(x_{-(n+1)})$
et son inverse $g$ est à valeurs dans $B(0,\rho)$ où $\rho:=\rho\big(x_{-(n+1)} \big)$ pour tout $\hat x$ 
d'un
ensemble $\widehat Y$ de $\hat\mu$-mesure pleine et tout entier $n$.
Considérons alors le développement de Taylor $g-d_0g=:\Sigma_{p\ge 2} Q_p$ où $Q_p$ désigne une application 
homogène
de degré $p$. Pour tout $u\in B(0,r)$ on a $\vert\vert Q_p(u)\vert\vert = \vert\vert\frac{1}{2\pi}
\int_{0}^{2\pi}g(e^{i\theta}u)e^{-ip\theta}d\theta \vert\vert \le \rho$ et donc

$$ \vert\vert (g-d_0 g) (u)\vert\vert \le \Sigma_{p\ge 2} \frac{\vert\vert u \vert\vert ^p}{r^p}
\vert\vert Q_p\big( \frac{ru}{\vert\vert u \vert\vert}\big)\vert\vert \le
 \rho\Sigma_{p\ge 2} \big(\frac{\vert\vert u \vert\vert }{r}\big)^p$$
 
 Si de plus $u \in d_0f^{-n}_{\hat x}\big[B(0,\gamma\eta({\hat x})\big]$ alors 
 $\frac{\vert\vert u \vert\vert}{r} \le \gamma e^{-n(\lambda_1 - \epsilon)}=: \delta_n$ 
 (cf. Lemme \ref{lempr}, 3.) 
 et donc
  $\vert\vert (g-d_0 g) (u)\vert\vert \le \frac{1}{1-\delta_1}\rho \delta_n^2$. La première assertion du 
  lemme \ref{lem3.4} s'en déduit car $\rho(x_{-(n+1)})$ a un taux de croissance exponentiel nul 
  (cf. Lemme
  \ref{lempr}, 2.). La seconde assertion découle presque immédiatement de la définition des exposants de Liapounov.
  \hfill$\square$\\

Nous chercherons, dans la section suivante, des conditions pour que les transformations linéaires 
$(d_0f^n_x)^{-1}$ intervenant dans la définition de la linéarisabilité en $x\in {\bf P}^k$ soient
assimilables à des homothéties. Cela revient à contr\^ oler le volume et la taille des ellipsoïdes
$(d_0f^n_x)^{-1}\big[B(0,1)\big]$. Le théorème de Briend-Duval, déjà implicitement  pour établir le lemme
\ref{lempr},  
majore le taux de décroissance exponentielle
de la taille 
de ces ellipsoïdes par $-\frac{Log d}{2}$. Cela signifie que, pour tout $\epsilon>0$, on a  
$\vert\vert (d_0f^n_x)^{-1} \vert\vert \lesssim e^{n\epsilon}d^{-\frac{n}{2}}$ pour $n$ assez grand.
En reprenant la méthode de Briend-Duval dans le
contexte de la proposition \ref{prop3.2}, nous allons établir une majoration plus précise : 
$\vert\vert (d_0f^n_x)^{-1} \vert\vert \lesssim d^{-\frac{n}{2}}$.
 A cet effet, nous introduisons pour tout $\tau>0$
les sous-ensembles de     
${\cal B}_n(\rho)$ suivants :

$${\cal L}{\cal B}_n(\rho,\tau):={\cal B}_n(\rho) \cap 
\{x\in {\bf P}^k / \vert\vert (d_0f^n_x)^{-1} \vert\vert
\le \tau d^{-\frac{n}{2}}\}$$          

Sous les hypothèses de la proposition \ref{prop3.2}, nous allons montrer que : 

$$\liminf_n \mu[{\cal L}{\cal B}_n(\rho,\tau)]\ge\alpha(\rho)-\frac{C}{\tau^2\rho^2}$$

où $C>0$ ne dépend que de $f$ et $\alpha: ]0,1] \to {\bf
R}^+$ est la fonction introduite à la proposition \ref{prop3.2}. Ceci établira le théorème \ref{THB}.\\

La preuve est basée sur le principe suivant. Par tout point de $\widehat{\cal L}{\cal B}_n(\rho,\tau)
:={\cal B}_n(\rho)\backslash{\cal L}{\cal B}_n(\rho,\tau)$ passe un disque dont le diamètre est
au moins égal à $\tau \rho d^{-\frac{n}{2}}$ et dont l'image par $f^n$ reste contenue dans une boule
de rayon fixé. Comme $f^{n*}T=d^n T$, il passe donc par tout point de $\widehat{\cal L}{\cal
B}_n(\rho,\tau)$ un grand disque peu chargé par $T$. Des techniques pluripotentialistes 
permettent alors de majorer précisément la masse de cet ensemble de points pour la mesure $\mu=T^k$.
Nous adoptons la définition suivante :\\

\noindent{\bf Définition} 
{\it On dit qu'un disque holomorphe $\sigma : \Delta 
\to {\bf C}^k$ est de taille $\alpha$ et passe par $z\in {\bf C}^k$ si il est 
de la forme $\sigma(u)=z+\alpha u. v
+\beta(u)$ où $\alpha>0$, $v$ est un vecteur unitaire de ${\bf C}^k$, $\beta(0)=0$
 et $\vert\vert \beta
\vert\vert\le \frac{\alpha}{1000}$  }.\\  

\underline{Preuve du théorème \ref{THB}} :\\

l'ingrédient principal est le théorème suivant dont la preuve sera résumée dans l'appendice.\\

\noindent{\bf Théorème} (Briend-Duval) 
{\it Soit $S:=dd^c w$ un $(1,1)$ courant positif fermé de potentiel $w$ continu sur $P(0,R)$ et 
$E\subset P(0,\frac{R}{2})$. On suppose que par tout $z \in E$ passe un disque holomorphe 
$\sigma_z : \Delta
\to {\bf C}^k$ de taille $\alpha$ et qu'il existe une fonction $h_z$ harmonique
telle que 
$\vert w \circ \sigma_z - h_z \vert \le \epsilon$ sur $\Delta$. Alors il existe une constante $C(w)$ ne
dépendant que de $w$ telle que $S^k(E) \le C(w) \frac{k^2}{\alpha^2}\epsilon$}.\\

En vue d'utiliser ce résultat, nous fixons des systèmes de coordonnées locales sur ${\bf P}^k$.
Considérons un recouvrement de ${\bf P}^k$ par des ouverts $U_1,...,U_N$ centrés en des points $m_j$
et tel que sur chaque $U_j$ nous puissions fixer des déterminations des cartes $\tau_x$ dépendant
différentiablement de $x$ (cf. Préliminaires).
Posons $\tau_j:=\tau_{m_j}$ puis, pour $R>0$ fixé, $V_j:=\tau_j(P(0,R))$. 
Si le recouvrement est assez
fin alors les propriétés suivantes sont satisfaites :

\begin{itemize} 

\item[(i)] $U_j \subset \tau_j(P(0,\frac{R}{2}))$ et $\tau_x(P(0,\frac{R}{2}))\subset V_j$ pour tout $x\in U_j$
\item[(ii)] $\forall x \in U_j, 
\vert\vert \tau_j^{-1}\circ \tau_x - (\tau^{-1}_j(x)+Id) \vert\vert_{{\cal C}^1,\overline{P(0,\frac{R}{2})}}
\le\frac{1}{1000}$ 
\item[(iii)]  
$\mu \{x \in U_j \cap {\cal B}_n(\rho)\;/\;     
(d_0f^n_x)^{-1}\big[B(0,\rho)\big] \subset P(0,\frac{R}{2})\}=\mu (U_j \cap {\cal B}_n(\rho)) - \epsilon_{n,j}$ avec 
$\lim_n \epsilon_{n,j}=0$

\end{itemize}

puis, si $R_0$ est pris assez petit ($R_0$ a été introduit au lemme \ref{lempr}) :

\begin{itemize}
\item[(iv)] $\forall x \in {\bf P}^k, \exists l \in \{1,...,N\}$ tel que $\tau_x\big[B(0,R_0)\big]
\subset V_l$ 
\end{itemize}

enfin, si $v_j$ désigne un potentiel continu de $T$ sur $V_j$

\begin{itemize}
\item[(v)]  $T = dd^c v_j$ et $\vert v_j\vert \le M$ sur $V_j$ pour tout $j \in \{1,...,N\}$. 
\end{itemize}

D'après la proposition \ref{prop3.2} on a $\mu({\cal B}_n(\rho))\ge\alpha(\rho)$ et, comme il s'agit 
de minorer
$\liminf_n\mu[{\cal L}{\cal B}_n(\rho,\tau)]$, la propriété (iii) montre que l'on peut considérer que :
\begin{equation}\label{inclu}
(d_0f^n_x)^{-1}\big[B(0,\rho)\big] \subset P(0,\frac{R}{2}), \;\;\forall x \in U_j \cap {\cal B}_n(\rho)
\end{equation}
Rappelons que 
$\widehat{\cal L}{\cal B}_n(\rho,\tau)
:={\cal B}_n(\rho)\backslash{\cal L}{\cal B}_n(\rho,\tau)$, pour tout $j \in \{1,...,N\}$
nous allons établir que : 
\begin{equation}\label{majo} \mu\big[{\widehat{\cal L}}{\cal B}_n(\rho,\tau)\cap U_j\big] \le
C(v_j\circ\tau_j)\frac{Mk^2}{\tau^2\rho^2}
\end{equation}
Soit donc 
$x\in{\widehat{\cal L}}{\cal B}_n(\rho,\tau)\cap U_j$ et $V_n(x)$ un vecteur unitaire réalisant 
$\vert\vert (d_0f^n_x)^{-1} \vert\vert$. On définit un disque affine 
$\Phi_{n,x}: \bar{\Delta}\to{\bf C}^k$
de diamètre au moins égal à $\tau \rho d^{-\frac{n}{2}}$ par  
$$\Phi_{n,x}(t):=(d_0f^n_x)^{-1}[t\rho. V_n(x)]=:t\rho. v_n(x).$$ Comme
$x\in U_j$, (\ref{inclu}) et (i) permettent
de définir un nouveau disque $\Phi_{j,n,x}:\Delta \to P(0,R)$ par     
$\Phi_{j,n,x}:=\tau_j^{-1}\circ \tau_x \circ \Phi_{n,x}$. Compte tenu de la propriété (ii),   
$\Phi_{j,n,x}$ est un disque holomorphe de taille $\alpha:=\rho\vert\vert v_n(x) \vert\vert 
\ge \tau \rho d^{-\frac{n}{2}}$
passant par $\tau_j^{-1}(x)$.\\
Choisissons $l \in \{1,...,N\}$ tel que $\tau_{f^n(x)}\big[B(0,R_0)\big]
\subset V_l$ (propriété (iv)) alors, comme $x \in {\cal B}_n(\rho)$, on a $f^n \circ \tau_x \circ 
\Phi_{n,x}(\Delta) \subset \tau_{f^n(x)}\big[B(0,R_0)\big]
\subset V_l$ et donc $$dd^c(v_l\circ f^n \circ \tau_x \circ 
\Phi_{n,x})= (f^n \circ \tau_x \circ 
\Phi_{n,x})^*T=(\tau_x \circ 
\Phi_{n,x})^*f^{n*}T=d^n(\tau_x \circ 
\Phi_{n,x})^*T.$$
 Par ailleurs, puisque $\tau_x\circ \Phi_{n,x}(\Delta)\subset V_j$ (cf. (\ref{inclu}) et (i)), on a         
$$(\tau_x \circ 
\Phi_{n,x})^*T=dd^c(v_j\circ \tau_x \circ 
\Phi_{n,x})=dd^c(v_j\circ \tau_j \circ 
\Phi_{j,n,x})$$
Ainsi, $dd^c\big[v_l\circ f^n \circ \tau_x \circ 
\Phi_{n,x} - d^n v_j\circ \tau_j \circ 
\Phi_{j,n,x}\big] =0$.
Autrement dit, la fonction entre crochets est harmonique sur $\Delta$ et, puisque $\vert v_l\vert\le M$, le
potentiel $w:=v_j\circ\tau_j$ de $\tau_j^* T =:S$ diffère d'au plus $\frac{M}{d^n}$ d'une fonction
harmonique $h_x$ sur le disque 
$\sigma_x:=\Phi_{j,n,x} $ de taille $\alpha\ge \tau \rho d^{-\frac{n}{2}}$. 
Dans ces conditions, (\ref{majo}) découle
immédiatement
du théorème de Briend-Duval. On en déduit l'estimation annoncée avec 
$C=Mk^2\Sigma_{1}^{N}C(v_j\circ\tau_j)$.  \hfill$\square$\\

Il est utile de disposer d'une version du théorème \ref{th3.1} où les orbites issues d'un borélien prescrit
sont assujetties à récurrence.
Cette précision s'obtient facilement gr\^ ace au caractère 
mélangeant de $\mu$ et conduit à l'énoncé du théorème \ref{THL}.\\

\underline{Preuve du théorème \ref{THL}} :\\

\noindent Posons ${\cal L}{\cal B}_n(\rho,\tau,B):={\cal L}{\cal B}_n(\rho,\tau)\cap B \cap f^{-n}(B)$
et ${\cal L}{\cal B}(\rho,\tau,B):=\limsup_n {\cal L}{\cal B}_n(\rho,\tau,B)$. Il est clair que si 
$x\in{\cal L}{\cal B}(\rho_0,\tau_0,B)$ alors il existe une suite extraite $(f^{n_j})_j$ vérifiant
les trois assertions du théorème \ref{THL}. Il suffit donc     
d'observer que $\mu \big({\cal L}{\cal B}_n(\rho,\tau,B) \big)\ge (1-0) \mu(B)^2$
pourvu que $\rho_0$,$\frac{1}{\tau_0}$ soient assez petits et $n$ assez grand.
Or ceci résulte
immédiatement du théorème \ref{THB} et du caractère mélangeant de $\mu$. \hfill$\square$

\section{Systèmes linéarisables.}

Dans cette partie nous commen\c cons par établir des conditions pour que la suite des itérées $(f^n)_n$ d'un 
système 
$({\bf P}^k, f, d, \mu)$ soit $\mu$-presque partout linéarisable par des homothéties de
rapport $d^{-\frac{n}{2}}$. Nous montrons ensuite que le courant de Green d'un tel système est une forme
lisse strictement positive sur un ouvert non vide de ${\bf P}^k$. D'après \cite{BL}, 
$f$ est alors un endomorphisme de Lattès.\\

Dans cette perspective, nous adoptons la définition suivante :\\

\begin{defi}\label{defi4.1}
On dit que le système $({\bf P}^k, f, d, \mu)$ est $\sqrt{d}$-linéarisable si pour $\mu$-presque tout
$x \in {\bf P}^k$ 
il existe $\nu(x)>0$ ainsi qu'une suite $[f^{n_j} \circ \tau_x \circ \big(d^{-{\frac{n_j}{2}}}Id\big)]_j$ 
qui converge
uniformément vers une limite injective sur $B(0,\nu(x))$.     
\end{defi} 

A la lumière du théorème \ref{THL}, ceci revient à exiger que les ellipsoïdes 
$(d_0f^{n_j}_x)^{-1}\big[B(0,1)\big]$ soient assimilables à des 
boules euclidiennes de rayon $d^{-\frac{n_j}{2}}$. Or, la taille de ces
ellipsoïdes étant au plus de l'ordre de $d^{-\frac{n_j}{2}}$, il s'agit en 
fait d'en contr\^ oler
le volume. Nous introduisons donc, pour tout $0<\nu<1$, les ensembles   
suivants où $J_0\;f^n_x$ désigne le Jacobien complexe de $f^n_x$ en $0$ :\\

$${\cal V}_n(\nu):=\{x\in {\bf P}^k  \;/\;
\nu^2d^{kn}\le \vert J_0\;f^n_x\vert^2 \le \frac{1}{\nu^2} d^{kn}\}$$

En estimant la masse des ${\cal V}_n(\nu)$, nous caractérisons la linéarisabilité d'un système 
$({\bf P}^k, f, d, \mu)$ par l'absence de partie singulière dans la mesure $\mu$ :

\begin{prop}\label{pr4.1}
Pour un système $({\bf P}^k, f, d, \mu)$ les propriétés suivantes sont équivalentes.
\begin{itemize} 

\item[1)] $\mu$ n'est pas singulière par rapport à la mesure de Lebesgue $m:=\omega^k$
\item[2)] $\mu$ est absolument continue par rapport à la mesure de Lebesgue $m:=\omega^k$
\item[3)]i) $\exists \beta : ]0,1] \to {\bf
R}^+$ telle que $\lim_{\nu \to 0} \beta(\nu) = 1$ et $\liminf_n\mu[{\cal V}_n(\nu)]\ge\beta(\nu)$
\item[\phantom{3)}]ii) les exposants du système sont tous égaux à $\frac{Log d}{2}$
\item[4)] Le système est $\sqrt{d}$-linéarisable.

\end{itemize}
\end{prop}

\underline{Preuve de la proposition \ref{pr4.1}} :\\

$1) \Rightarrow 2)$. Soit $\mu = \mu_a + \mu_s$ la décomposition de Lebesgue de $\mu$. On a $\mu = f_{*}\mu= 
f_{*}\mu_a + f_{*}\mu_s$. Or, $f$ étant un rev\^ etement ramifié, les images directes ou réciproques par $f$
d'ensembles
$m$-négligeables sont encore $m$-négligeables. Il résulte alors facilement des définitions que
$f_{*}\mu_a \ll m$ et $f_{*}\mu_s\bot m$. Donc, la décomposition de Lebesgue de $\mu$ étant unique, les mesures 
$\mu_a$ et $\mu_s$ sont invariantes : $f_{*}\mu_a = \mu_a$, $f_{*}\mu_s=\mu_s$. Etant ergodique, la mesure 
$\mu$ est extrémale parmi les mesures de probabilité invariantes par $f$. Donc, si $m_a:=\mu_a({\bf P}^k)>0$
alors $m_s:=\mu_s({\bf P}^k)=1- m_a$ et l'identité $\mu =m_a \frac{\mu_a}{m_a} +(1-m_a) \frac{\mu_s}{m_s}$
force l'égalité $\mu=\mu_a$.\\

$2) \Rightarrow 3)i)$. Puisque $\mu \ll m$, il existe une fonction $\varphi$ $m$-intégrable sur ${\bf P}^k$ telle que 
$\mu=\varphi
m$. De plus, par le théorème de Lusin, on trouve pour tout $n\in {\bf N}$ des
fonctions continues $g_n$,$h_n$ et des boréliens $C_n(\varphi)$,$C_n(\varphi \circ f^n)$ tels que : 
$$ \varphi = g_n \;sur\; C_n(\varphi) \;et\; \mu\big[C_n(\varphi)\big] \ge 1-\frac{1}{n}$$    
$$ \varphi \circ f^n = h_n \;sur\; C_n(\varphi \circ f^n) \;et\; \mu\big[C_n(\varphi \circ f^n)\big] 
\ge 1-\frac{1}{n}\;.$$
Pour tout $0<\nu <1$ définissons un borélien $A_{\nu}$ par 
$A_{\nu}:=\{x \in {\bf P}^k\; /\; \nu <\varphi <\frac{1}{\nu}\}$
et posons $\beta(\nu):=\mu(A_{\nu})^2$. 
Soit alors $Z_{n,\nu}:=\big[f^{-n}(A_{\nu})\cap A_{\nu}\big] \cap \big[C_n(\varphi) \cap 
C_n(\varphi \circ f^n)\big]
\cap Y^{c}$ où $Y:=\cup_p Crit\;f^p$ et $Z_{n,\nu}^{Leb}$ l'ensemble des points de Lebesgue de 
$Z_{n,\nu}$ c'est à dire :
$$Z_{n,\nu}^{Leb}:=\{x\in Z_{n,\nu}\;/\;\lim_{s\to 0}
\frac{m\big[B(x,s)\cap Z_{n,\nu}\big]}{m\big[B(x,s)\big]}=1 \}$$
Puisque $\mu\ll m$ on a $\mu\big(Z_{n,\nu}^{Leb}\big)=\mu\big(Z_{n,\nu}\big)$. Alors, compte tenu du caractère
mélangeant de $\mu$ et du fait que $\mu(Y)=0$, on voit que pour $n$ assez grand on a :
$$\mu \big(Z_{n,\nu}^{Leb}\big)\ge \big[\beta(\nu)-0\big] - \frac{2}{n}$$
Il suffit donc de vérifier que $Z_{n,\nu}^{Leb}\subset {\cal V}_n(\nu)$. Fixons $x \in Z_{n,\nu}^{Leb}$. 
Puisque $ x \notin Crit\;f^n$, il existe $s_0>0$ tel que $f^n$ soit injective sur $B(x,s_0)$. En outre, $x$
étant un point de Lebesgue de $Z_{n,\nu}$, on peut diminuer $s_0$ de fa\c con à ce que       
$m\big[B(x,s)\cap Z_{n,\nu}\big] \ge \frac{1}{ 2}m\big[B(x,s)\big]>0$ pour tout $0<s<s_0$. Par changement de
variables, d'abord par rapport à $\mu=\varphi m$ qui est de Jacobien constant égal à $d^k$ puis par rapport à
$m=\omega^k$, on obtient :
$$d^{kn}\int_{B(x,s)\cap Z_{n,\nu}}\varphi m = \int_
{f^n[B(x,s)\cap Z_{n,\nu}]}\varphi m
=\int_{B(x,s)\cap Z_{n,\nu}}\varphi \circ f^n (f^{n*}\omega^k)$$
Or, puisque $C_n(\varphi) \cap 
C_n(\varphi \circ f^n)$ contient $Z_{n,\nu}$, on peut remplacer $\varphi$ par $g_n$ et $\varphi \circ f^n$ par $h_n$ dans ces
intégrales. Après normalisation par $m(s,n,\nu):= m[B(x,s)\cap Z_{n,\nu}]$ cela donne :
$$\frac{d^{kn}}{m(s,n,\nu)}\int_{B(x,s)\cap Z_{n,\nu}}g_n m =
\frac{1}{m(s,n,\nu)} 
\int_{B(x,s)\cap Z_{n,\nu}}h_n (f^{n*}\omega^k)$$ 
Faisons tendre $s$ vers $0$, comme les fonctions $g_n$,$h_n$ sont continues et 
$(f^{n*}\omega^k)_x=\vert J_0\;f^n_x\vert^2 (\omega^k)_x$, on obtient :
$$d^{kn} \varphi(x) = d^{kn} g_n(x)=
h_n(x)\vert J_0\;f^n_x\vert^2 = \varphi \circ f^n(x)\vert J_0\;f^n_x\vert^2 $$
c'est à dire $\frac{\vert J_0\;f^n_x\vert^2}{d^{kn}}=
\frac{\varphi(x)}{\varphi \circ f^n(x)}$ ce qui, puisque $x$ et $f^n(x)$ appartiennent à $A_{\nu}$, 
montre que $x \in {\cal V}_n(\nu)$.\\

$3)i) \Rightarrow 3)ii)$. On sait que $\lim_n \frac{1}{n}Log \vert J_0\;f^n_x\vert^2 =2\sum_{i=1}^{k}
\lambda_i$ pour $\mu$-presque tout $x$. Par ailleurs, si $x\in {\cal V}(\nu):=\limsup_n {\cal V}_n(\nu)$
on a $\lim_n \frac{1}{n}Log \vert J_0\;f^n_x\vert^2 = k Log d$. Or, d'après 3)i), $\mu[{\cal V}(\nu)]
\ge \beta(\nu) \ge \frac{1}{2}$ pourvu que $\nu$ soit assez petit. On a donc   
$\sum_{i=1}^{k} \lambda_i = k \frac{log d}{2}$, d'où 3)ii) puisque $\lambda_i \ge \frac{log d}{2}$
.\\

$3) \Rightarrow 4)$. On peut, compte tenu de 3)ii), appliquer le théorème \ref{THB}. Nous en reprenons les
notations et posons ${\cal LVB}_n(\rho,\tau,\nu):= {\cal LB}_n(\rho,\tau)\cap 
{\cal V}_n(\nu)$, ${\cal LVB}(\rho,\tau,\nu):= \limsup_n {\cal LVB}_n(\rho,\tau,\nu)$.
D'après 3)i) et le théorème \ref{THB}, 
$\mu[{\cal LVB}(\rho,\tau,\nu)]$ est arbitrairement proche de $1$ pourvu que $\rho$,$\nu$ soient assez
petits et $\tau$ assez grand. Il suffit donc de montrer que $(f^n)_n$ est linéarisable par $\Lambda_n:=
(\sqrt{d})^{-n} Id$ lorsque $x \in {\cal LVB}(\rho,\tau,\nu)$. Soit donc $(n_j)_j$ une suite strictement
croissante d'entiers telle que $x\in{\cal LVB}_{n_j}(\rho,\tau,\nu)$ pour tout $j$. Puisque     
${\cal LVB}(\rho,\tau,\nu)  \subset {\cal B}(\rho)$ on a 
$f^{n_j}_x \circ (d_0f^{n_j}_x)^{-1} \big(B(0,\rho)\big)\subset
B(0,R_0)$ pour tout $j$ et il nous reste donc à voir que $(d_0f^{n_j}_x)^{-1}$ est équivalente à 
$\Lambda_{n_j}$. A cet effet, notons $\delta_{j,1}\le ...\le \delta_{j,k}$ les valeurs singulières de  
$(d_0f^{n_j}_x)^{-1}$. On a $\delta_{j,k}\le \tau (\sqrt{d})^{-n_j}$ car $x \in {\cal LB}_{n_j}(\rho,\tau)$
et $(\delta_{j,1}...\delta_{j,k})^2 = \vert J_0\;f^{n_j}_x \vert^{-2} \ge \nu^2 (d)^{-kn_j}$ car
$x\in {\cal V}_{n_j}(\nu)$. On
en déduit l'équivalence voulue : $\nu (\tau)^{1-k} (\sqrt{d})^{-n_j} \le\delta_{j,1}\le ...\le \delta_{j,k}\le 
\tau (\sqrt{d})^{-n_j}$.\\ 

$4) \Rightarrow 1)$. Soit $x\in {\bf P}^k$ un point où 
$(f^n)_n$ est linéarisable par $\Lambda_n:=
(\sqrt{d})^{-n} Id$. Cela signifie qu'il existe $\rho>0$ et une suite strictement croissante d'entiers  
$(n_j)_j$ tels que $f^{n_j}_x \circ \Lambda_{n_j} : B(0,\rho)\to
B(0,R_0)$ soit une suite d'injections. Notons $B_{n_j}:=B(0,\rho d^{\frac{-n_j}{2}})$ alors, puisque
$f^*\mu=d^k\mu$ on a $\mu [f^{n_j} \circ \tau_x (B_{n_j})] = d^{kn_j} \mu [\tau_x (B_{n_j})]$. Il s'ensuit que
$$\liminf_{r\to 0} \frac{\mu [\tau_x (B(0,r))]}{m [\tau_x (B(0,r))]}\le
\liminf_j \frac{\mu [\tau_x (B_{n_j})]}{m [\tau_x (B_{n_j})]}\lesssim 
\liminf_j \frac{\mu [\tau_x (B_{n_j})]}{d^{-kn_j}}$$
$$=\liminf_j \mu [f^{n_j} \circ \tau_x (B_{n_j})] \le 1$$

et donc $\mu \ll m$, puisque ceci est vrai pour $\mu$-presque tout $x$. \hfill$\square$    

\begin{rem}\label{remrec}
Comme nous l'avons fait pour établir le théorème \ref{THL}, une légère modification de l'argumentation
dans $3) \Rightarrow 4)$ permet de choisir la sous-suite 
$[f^{n_j} \circ \tau_x \circ \big(d^{-{\frac{n_j}{2}}}Id\big)]_j$ de 
fa\c con à ce que $f^{n_j}(x)$
ne s'échappe pas d'un borélien de $\mu$-mesure strictement positive prescrit.  
\end{rem}

La fin de cette section est dévolue à la preuve du théorème \ref{TH} proprement dite. Au vu de la 
proposition \ref{pr4.1},
il s'agit de caractériser les systèmes $({\bf P}^k, f, d, \mu)$ qui sont  
$\sqrt{d}$-linéarisables. Comme le courant de Green $T$ satisfait les équations fonctionnelles
 $f^{n*} T=d^nT$, nous
verrons que pour de tels systèmes le procédé de linéarisation peut fournir des coordonnées locales dans
lesquelles
$T$ est une forme lisse définie positive. L'endomorphisme $f$
est alors un exemple de Lattès en vertu du résultat de \cite{BL}.\\

\underline{Preuve du théorème \ref{TH}} :\\

Le procédé de linéarisation permet d'établir le lemme suivant :\\

\begin{lem}\label{lem4.1}
Soit $({\bf P}^k, f, d, \mu)$ un système $\sqrt{d}$-linéarisable et 
$S$ un $(1,1)$ courant positif sur ${\bf P}^k$ tel que $f^{*}S=dS$ ($S$ n'est pas nécessairement fermé).\\
1) Si $S=S_a$ sur un ouvert $\Omega \subset {\bf P}^k$ chargé par $\mu$ alors il existe 
une boule $B(0,r)\subset {\bf C}^k$ et un biholomorphisme
$\Phi:B(0,r) \to \Omega' \subset \Omega$ tels que $\Phi^{*}S$ soit une forme différentielle à coefficients 
constants et $\mu(\Omega')>0$.\\   
2) Soit  
$\Omega$ un ouvert de ${\bf P}^k$ chargé par $\mu$ et sur lequel $S$ dérive d'un potentiel $p.s.h$ continu
$v$ $(S=dd^c v)$. Si $S_a$ est identiquement nul sur $\Omega$ alors $\mu(\Omega \cap\; Supp\;S) =0$\\ 
\end{lem}

Soit $\Omega$ un ouvert de ${\bf P}^k$ tel que $\mu(\Omega)>0$. Considérons la décomposition $T=T_a+T_s$.
La première assertion du lemme \ref{lem4.1} appliquée à 
$T_a$ permet de supposer que $T_a\vert_{\Omega}$ est donné par une forme
$H$ à coefficients constants dans de bonnes coordonnées. En particulier cela montre que $T_a$ possède un
potentiel continu sur $\Omega$ et il en va donc de m\^ eme pour $T_s=T-T_a$. Ceci permet, sur 
$\Omega$, d'exprimer
$\mu$ sous la forme d'une somme de mesures positives obtenues comme
 produits extérieurs de $T_a$ et $T_s$ :
\begin{equation}\label{decmu}
\mu=T^k=\big(T_a+T_s\big)^k=T_a^k + \sum_{j=1}^k C_k^j\;\; T_s^j \wedge T_a^{k-j}     
\end{equation}
Puisque $(T_s)_a \equiv 0$, la seconde assertion du lemme \ref{lem4.1} montre que $\mu$ ne charge pas 
$\Omega\cap Supp\;T_s$ et
donc, au vu de (\ref{decmu}), la mesure $T_a^k$ n'est pas identiquement nulle sur $\Omega$. 
Autrement dit la forme $H$ n'est pas
dégénérée. Par ailleurs, chaque terme du second membre de (\ref{decmu}) doit, en tant que mesure positive, 
\^ etre
absolument continue. En particulier, $\big(T_s \wedge T_a^{k-1}\big)_s$ est identiquement nul sur $\Omega$. 
Or, puisque
$H$ est strictement positive, $\big(T_s \wedge T_a^{k-1}\big)$ est équivalente à la mesure trace
$\sigma_{T_s}$ de $T_s$ et donc $\sigma_{T_s}$ est nulle sur $\Omega$. Ainsi,
$T_s\vert_{\Omega}\equiv 0$ et $T$ coïncide sur $\Omega$ avec une forme lisse définie positive.
\hfill$\square$\\

\underline{Preuve du lemme \ref{lem4.1}} :\\

1) Quitte à diminuer $\Omega$ on peut supposer que  
$\tau_x^{-1}$ soit défini sur $\Omega$ pour tout $x\in\Omega$. Choisissons $x_0\in\Omega$. Alors,
puisque $S\equiv S_a$, on peut écrire
$\tau_{x_0}^{*}S$ sous la forme
$$\tau_{x_0}^{*}S=\frac{i}{2}\sum_{1\le p,q\le k} h_{p,q}(z)\;dz_p\wedge d\bar z_q$$
où les $h_{p,q}$ sont des fonctions $L^1_{loc}$.
Soit $\cal C$ l'ensemble des points de $\tau_{x_0}^{-1}(\Omega)$ où chacune des fonctions $h_{p,q}$
est continue en mesure et $\cal R$ l'ensemble des points de $(\Omega\cap Supp\;\mu)$ où $(f^n)_n$ est
linéarisable
par des homothéties $\Lambda_n:=d^{-\frac{n}{2}} Id$. Comme $\mu\ll m$, on voit gr\^ ace à la proposition
 \ref{pr4.1} 
que $\mu\big[\tau_{x_0}({\cal C})\cap{\cal R}\big]>0$. Soit alors 
$x_1 \in \tau_{x_0}({\cal C})\cap{\cal R}$ et $\psi:=\tau_{x_0}^{-1}
\circ \tau_{x_1}$. Comme $\tau_{x_1}^{*}S = \psi^{*} \tau_{x_0}^{*}S=
\frac{i}{2}\sum_{1\le p,q\le k} h_{p,q}\circ \psi \;d\psi_p\wedge d\bar {\psi}_q$ on voit que, quitte à 
remplacer $x_0$ par $x_1$, on peut supposer que $0\in {\cal C}$ et $x_0 \in {\cal R}$.\\ 
Soit donc $\Phi_n:=f^n \circ \tau_{x_0} \circ
\Lambda_n$. Modulo extraction, $\Phi_n$
converge vers un biholomorphisme $\Phi : B(0,\nu)\to\Omega'$ et l'on peut supposer que $\Phi_n(0)=f^n(x_0)$
reste dans $V_{x_0}\cap Supp\;\mu$ où $V_{x_0}$ est un voisinage arbitrairement petit de $x_0$ (cf. Remarque
\ref{remrec}). Ainsi $\Phi(0)\in \Omega \cap Supp\;\mu$ et, quitte à diminuer $\nu$, on a $\Omega'\subset
\Omega$ et $\mu(\Omega')>0$.
D'après l'invariance de $S$ on a :    
$$\Phi_n^{*}S=\Lambda_n^{*} \tau_{x_0}^{*}f^{n*}S=d^{n} \Lambda_n^{*} [\tau_{x_0}^{*}S]=
\frac{i}{2}\sum_{1\le p,q\le k} h_{p,q}\circ \Lambda_n\; dz_p\wedge d\bar z_q$$
d'où, en passant à la limite, $\Phi^{*}S=\frac{i}{2}\sum_{1\le p,q\le k} h_{p,q}(0)\;dz_p\wedge d\bar z_q
=:H$.\\

2) Quitte à diminuer $\Omega$ on peut supposer que $S=dd^c v$ sur un voisinage $V$ de $\overline{\Omega}$. Supposons
que $\mu(\Omega \cap Supp\;S) >0$. Soit $\Lambda_n$ l'homothétie de rapport $d^{-\frac{n}{2}}$. D'après la
proposition 
\ref{pr4.1} et la remarque \ref{remrec}, il existe ${\cal R} \subset (\Omega \cap\; Supp\;S)$ tel que   
$\mu({\cal R}) >0$ et, pour tout point $x_0\in {\cal R}$, il existe une suite  
$\Phi_{n_j}:=f^{n_j} \circ \tau_{x_0} \circ
\Lambda_{n_j}$ qui converge uniformément sur $B(0,\nu(x_0))$ vers un biholomorphisme $\Phi$
telle que $f^{n_j}(x_0) \in {\cal R}$ pour tout $j$.

Nous allons montrer que $(\sigma_S)$ possède une dérivée de Radon-Nykodym strictement positive en
tout point de $\cal R$. Comme $\mu({\cal R}) >0$ et $\mu\ll m$, 
cela montrera que $\sigma_{S_a}$ (qui est égale à $(\sigma_S)_a$) et donc $S_a$ ne sont pas nuls sur
$\Omega$.\\
Comme $\Phi(0) \in \overline{\cal R} \subset \overline{\Omega}$ on peut diminuer
$\nu$ de fa\c con à ce que $\Phi_{n_j}\big(B(0,\nu)\big)$ et $\Phi\big(B(0,\nu)\big)$ soient contenus dans $V$.    
Soit $S_0:=\tau_{x_0}^{*} S$ et $\omega_0:=\frac{i}{2} dd^c\vert\vert z\vert\vert^2$. Il suffit de contr\^ oler 
la dérivée
de $\big(\sigma_{S_0}\big)_a$ à l'origine, c'est à dire d'établir que 
$$\limsup_n d^{kn}
\int_{B\big(0,d^{-\frac{n}{2}}\nu\big)}S_0\wedge\omega_0^{k-1} >0.$$ 
Or, puisque par hypothèse $f^{*}S=dS$, 
il vient

$$d^{kn_j}\int_{B(0,d^{-\frac{n_j}{2}}\nu)}S_0\wedge\omega_0^{k-1}  =
d^{(k-1)n_j}\int_{\Lambda_{n_j} [B(0,\nu)]}\tau_{x_0}^{*}f^{{n_j}*}S\wedge\omega_0^{k-1}=$$
$$d^{(k-1)n_j}\int_{B(0,\nu)}\Phi_{n_j}^{*}S\wedge\big(\Lambda_{n_j}^{*}\omega_0\big)^{k-1}=
\int_{B(0,\nu)}\Phi_{n_j}^{*}S\wedge\omega_0^{k-1}=
\int_{B(0,\nu)}dd^c(v\circ\Phi_{n_j})\wedge\omega_0^{k-1}$$

d'où par le théorème de convergence dominée :
$$\limsup_j d^{kn_j}\int_{B(0,d^{-\frac{n_j}{2}}\nu)}S_0\wedge\omega_0^{k-1} \ge
\int_{B(0,\nu)}dd^c(v\circ\Phi)\wedge\omega_0^{k-1}=
\int_{B(0,\nu)}\Phi^{*}S\wedge\omega_0^{k-1}.$$
 
\noindent Ceci achève la preuve car $\Phi(0) \in \overline{\cal R}\subset Supp\;S$ entra\^ \i ne 
$\int_{B(0,\nu)}\Phi^{*}S\wedge\omega_0^{k-1}=
\sigma_{\Phi^* S}\big[B(0,\nu)\big]>0$.\hfill$\square$\\
\bigskip
\bigskip
\bigskip

\noindent{\bf APPENDICE}\\

$\bullet$ Pour la commodité du lecteur, nous résumons la preuve du théorème de Briend-Duval que 
nous avons utilisé
dans la section 3.\\

\noindent{\bf Théorème} (Briend-Duval) 
{\it Soit $S:=dd^c w$ un $(1,1)$ courant positif fermé de potentiel $w$ continu sur $P(0,R)$ et 
$E\subset P(0,\frac{R}{2})$. On suppose que par tout $z \in E$ passe un disque holomorphe 
$\sigma_z : \Delta
\to {\bf C}^k$ de taille $\alpha>0$, c'est à dire de la forme $\sigma(u)=z+\alpha u. v
+\beta(u)$ où $v$ est un vecteur unitaire de ${\bf C}^k$, $\beta(0)=0$
 et $\vert\vert \beta
\vert\vert\le \frac{\alpha}{1000}$, tel qu'il existe une fonction $h_z$ harmonique
vérifiant 
$\vert w \circ \sigma_z - h_z \vert \le \epsilon$ sur $\Delta$. Alors il existe une constante $C(w)$ ne
dépendant que de $w$ telle que $S^k(E) \le C(w) \frac{k^2}{\alpha^2}\epsilon$}.\\

Soit $p_l$ la projection sur le $l$-ième axe de ${\bf C}^k$ et
$E_l:=\{z\in E / \vert\vert p_l(v_z) \vert\vert \ge \frac{1}{\sqrt{k}}\}$, on a $E=\cup_{l=1,k}E_l$. 
Pour fixer les idées estimons
$S^k(E_1)$.
On recouvre le polydisque $P(0,\frac{1}{2}R)$ par environ $N:=\frac{1}{4}
\frac{100 k}{\alpha^2}$ ellipsoïdes contenus dans $P(0,R)$ et de la forme $z+{\cal D}\big[
B(0,R)\big]$ où ${\cal D}(z_1,z')=\big(\frac{\alpha}{10\sqrt{k}}z_1,z'\big)$.\\

Soit ${\cal E}$ l'un de ces ellipsoïdes. Puisque ${\cal E}$ est strictement pseudoconvexe, il existe une
fonction $\hat w$ $p.s.h$ maximale sur ${\cal E}$, continue sur $\overline {{\cal E}}$ et coïncidant avec
$w$ sur $b{\cal E}$.\\
Si $z\in {\cal E}\cap E_1$, on voit facilement que le disque $\sigma_z(\Delta)$ traverse ${\cal E}$ au sens
où la composante connexe de $\sigma_z^{-1}\big({\cal E}\cap\sigma_z(\Delta)\big)$ passant par
l'origine ($C_0$)
est relativement compacte dans $\Delta$. Un argument de principe du maximum montre que $C_0$
est simplement connexe. En exhaustant $C_0$ par des domaines à bord suffisamment régulier on peut paramétrer
des disques holomorphes contenus dans $\cal E$ et dont le bord est arbitrairement proche de $b{\cal E}$.
Plus précisément, $\epsilon >0$ étant fixé, on trouve une
transformation conforme et continue au bord $\psi : \overline{\Delta}\to \psi(\overline{\Delta})\subset{\cal
E}$ telle
que $\psi(0)=0$ et
$\vert \hat w - w \vert \le \epsilon$ sur $\sigma_z\circ\psi(b\Delta)$. Posons $\tilde{\sigma}_z:=\sigma_z \circ
\psi$, soit $\tilde h$ la fonction harmonique sur $\Delta$ continue sur $\overline{\Delta}$ et
coïncidant avec $w\circ \tilde{\sigma}_z$ sur $b\Delta$. On a alors :

\begin{equation}\label{enca}
w(z)\le \hat w (z)\le \tilde h (0) +\epsilon
\end{equation}
la première inégalité découlant de la maximalité de $\hat w$ sur $\cal E$ et la seconde du principe du
maximum appliqué à $\hat w\circ\tilde{\sigma}_z - \tilde  h$ qui coïncide avec 
$\hat w\circ\tilde{\sigma}_z - w\circ\tilde{\sigma}_z$ sur $b\Delta$.\\

Par hypothèse on a $h_z\circ\psi -\epsilon \le 
w\circ\sigma_z\circ \psi=w\circ\tilde{\sigma}_z
\le h_z\circ\psi +\epsilon$ sur $\overline \Delta$. 
On a donc aussi $h_z\circ\psi -\epsilon \le \tilde h \le h_z\circ\psi +\epsilon$ et il s'ensuit que :
\begin{equation}\label{maj2}
\vert w(z) - \tilde h(0)\vert \le 2\epsilon.
\end{equation}
Les inégalités (\ref{enca}) et (\ref{maj2}) montrent 
que :
$${\cal E}\cap E_1 \subset {\cal E}(w,\epsilon):=\{
z\in {\cal E} / 0\le \hat w(z) -w(z) \le 3\epsilon \}.$$ La majoration annoncée résulte alors 
immédiatement de
l'estimation suivante qui est au coeur de la démonstration de Briend-Duval et pour laquelle nous renvoyons à 
\cite{BD} ou \cite{S} page 180, Théorème A.10.2 :\\

{\it Il existe une constante} $C(w)>0$  {\it telle que} $(dd^c w)^k
\big[{\cal E}(w,\epsilon)\big]\le C(w) \epsilon$.
\hfill$\square$\\
\\

$\bullet$ Nous reprenons, \emph{mutatis mutandis}, des arguments développés par Binder et DeMarco \cite{BdM}
dans le cas d'endomorphismes polynomiaux de ${\bf C}^k$ pour justifier le résultat suivant :\\

\noindent{\bf Théorème}  
{\it Soit $({\bf P}^k, f, d, \mu)$ un système dynamique, d'exposants 
$\lambda_1 \leq \cdots \leq \lambda_k$. La dimension de Hausdorff de $\mu$ vérifie : 
$Dim_H(\mu) \leq 2(k-1) + {\log d \over \lambda_k}$ }.\\

Rappelons que la dimension est définie comme la borne inférieure des dimensions de Hausdorff des 
boréliens de mesure totale. Ce résultat montre que si la dimension de $\mu$ est égale à $2k$, alors 
tous les exposants de $\mu$ sont minimaux, égaux à $\log d/2$.\\

Esquissons maintenant la démonstration : il s'agit d'exhiber 
pour tout $\epsilon > 0$ un borélien $Y$ de mesure totale vérifiant :
\begin{equation}\label{dim1}
\dim_H (Y) \leq 2(k-1)  + {\log d \over \lambda_k} +  {2k \over \lambda_k}.\epsilon 
\end{equation}

Soit  $\widehat A$ l'ensemble des points 
$\hat x = (x_{-n})_{n \geq 0}$ de l'extension naturelle $\widehat \PP$ vérifiant pour tout $n \geq 0$ : 
\[ f^{-n}_{\hat x} [ B(x_0 ,r_0) ] 
\supset B(x_{-n} , {r_0 \over \kappa_0}.e^{-n(\lambda_k + \epsilon)} ) 
\ \textrm{ et } \ Vol \left(  f^{-n}_{\hat x} [ B(x_0 ,r_0) ] \right) 
\leq \kappa_0 . e^{-2n(\lambda_1+ \cdots + \lambda_k) + n \epsilon} \]
On montre que 
$\hat \mu (\widehat A)>0$ pour $\kappa_0$ assez grand et $r_0$ assez petit 
(cf. \cite{BdM}, lemme 2).\\
 Soit $\widehat {A_n} := \hat f ^{-n} \widehat A$. 
La mesure $\hat \mu$ étant ergodique, le théorème de Birkhoff entra\^ \i ne que 
$\widehat Y := \limsup_n \widehat {A_n}$ est de mesure totale. 

On pose alors $Y := \pi_0 (\widehat Y)$ et $A_n := \pi_0 (\widehat {A_n})$, de sorte que 
$Y$ est aussi de mesure totale et est contenu dans $\limsup_n {A_n}$. Estimer la dimension de Hausdorff de $Y$ revient à estimer celle des ensembles $A_n$, pour $n$ assez grand. Par définition de $\widehat A$, tout point $y$ de $A_n$ vérifie :
\begin{enumerate}
\item
$f^n$ admet une branche inverse $g$ sur $B(f^n (y) ,r_0)$, telle que $g (f^n(y)) = y$.

\item
${\cal P} := g [ B(f^n (y) ,r_0) ] \supset B(y , {r_0 \over \kappa_0}.e^{-n(\lambda_k + \epsilon)} ) $

\item
$Vol ({\cal P} )  \leq k_0 . e^{-2n(\lambda_1+ \cdots + \lambda_k) + n \epsilon}$
 
\end{enumerate}

Il découle de ces propriétés que $A_n$ est recouvert par une famille 
$({\cal P}_i)_{i \in I}$ d'ouverts du type ${\cal A}$ dont le cardinal est de l'ordre de
$d^{kn}$. Il suffit pour le voir de recouvrir $\overline{A_0}$ par un nombre fini de boules
$B(x_{i_0},\frac{1}{4}r_0)$ puis d'observer que tout $y\in A_n$ est dans 
$g\big[B(x_{i_0},\frac{1}{2}r_0)\big]$ dès lors que $f^n(y)\in B(x_{i_0},\frac{1}{4}r_0)$. 
D'après le point 3, 
le volume de la réunion des ${\cal P}_i$ n'excède pas 
$d^{kn}. e^{-2n(\lambda_1+ \cdots + \lambda_k) + n \epsilon}$. 

Considérons à présent un recouvrement $({\cal M}_j)_{j \in J}$ de $A_n$ par des sous-ensembles de 
diamètre ${r_0 \over 100  \kappa_0}.e^{-n(\lambda_k + \epsilon)}$ provenant d'un maillage de ${\bf P}^k$. 
D'après le point 2, 
un sous-ensemble ${\cal M}_j$ intersectant $A_n$ est nécessairement contenu dans 
$\cup_{i \in I} {\cal P}_i$. On a donc :   
\[ Card(J) \leq {   Vol ( \cup_{i \in I} {\cal P}_i )  \over Vol({\cal M}_j) } 
\lesssim { d^{kn}. e^{-2n(\lambda_1+ \cdots + \lambda_k) + n \epsilon}   
\over \left(e^{-n(\lambda_k + \epsilon)} \right)^{2k} } \]
En minorant les exposants $\lambda_1, \cdots, \lambda_{k-1}$ par $\log d/2$, 
on obtient :
\[ Card(J) \lesssim  d^n . e^{n \left(  [2(k-1)\lambda_k] + (2k+1)\epsilon \right) }  \]
Il s'ensuit que la mesure de Hausdorff de $A_n$ de dimension 
$\alpha_\epsilon = 2(k-1) + \log d / \lambda_k + 2k \epsilon /\lambda_k$ est minorée $e^{-n \epsilon}$, 
pour $n$ assez grand. La  $\alpha_\epsilon$-mesure de Hausdorff de $Y \subset \limsup_n A_n$ est 
donc finie pour tout $\epsilon >0$, ce qui prouve le théorème.

\bigskip
\bigskip
\bigskip

{\footnotesize F. Berteloot}\\
{\footnotesize Universit\'e P. Sabatier, Toulouse III}\\
{\footnotesize Lab. Emile Picard, Bat. 1R2, UMR 5580}\\
{\footnotesize 118, route de Narbonne }\\
{\footnotesize 31062 Toulouse Cedex France}\\
{\footnotesize berteloo@picard.ups-tlse.fr}\\

\bigskip
\bigskip

{\footnotesize C. Dupont}\\
{\footnotesize Universit\'e Paris-Sud}\\
{\footnotesize Mathématique, Bat. 425, UMR 8628}\\
{\footnotesize 91405 Orsay, France}\\
{\footnotesize christophe.dupont@math.u-psud.fr}\\


\begin{thebibliography}{1}





\bibitem{BL} F. {Berteloot}, J.J. {Loeb}, \it{Une caractérisation géométrique des exemples de Lattès 
de ${\bf P}^k$},  \rm Bull. Soc. Math. Fr., {\bf 129}  \rm (2001), no. 2, 175-188. 

\bibitem{BdM} I. Binder, L. DeMarco, \it{Dimension of pluriharmonic measure and polynomial 
endomorphisms of ${\bf C}^n$}, \rm Int. Math. res. Not., {\bf 11} \rm(2003), 613-625. 

\bibitem{BD} J.Y. Briend, J. Duval, {\it Exposants de Liapounoff et distribution des points 
périodiques d'un endomorphisme de} ${\bf P}^k$, \rm Acta Math., \bf{182} \rm (1999), no. 2, 143-157. 

\bibitem{BD2} J.Y. Briend, J. Duval, \it{Deux caractérisations de la mesure d'équilibre d'un 
endomorphisme de ${\bf P}^k$}, \rm Publ. Math. Inst. Hautes Études Sci., \bf{93} \rm (2001), 145-159. 

\bibitem{DS} T.C. Dinh, N. Sibony, \it{Sur les endomorphismes holomorphes permutables de ${\bf P}^k$}, 
\rm Math. Ann., \bf{324} \rm (2002), no. 1, 33-70.

\bibitem{Du} C. Dupont, \it {Propriétés extrémales et caractéristiques des exemples de Lattès}, 
\rm Thèse de doctorat de l'Université Paul Sabatier, Toulouse, (2002).

\bibitem{Du1} C. Dupont, \it {Exemples de Lattès et domaines faiblement sphériques.}, \rm Manuscripta
Mathematica, (to appear). 

\bibitem{D} C. Dupont, \it {Endomorphismes holomorphes de ${\bf P}^k$ vérifiant la formule de Pesin}, 
\rm preprint.



\bibitem{FS} J.E. Fornaess, N. Sibony, \it{Complex Dynamics in higher dimensions}, 
\rm in Complex potential theory (Montréal, PQ, 1993), NATO ASI series Math. and Phys. Sci., \bf{439}
\rm, Kluwer Acad. Publ. (1994), 131-186.

\bibitem{FS2} J.E. Fornaess, N. Sibony, \it{Complex Dynamics in higher dimensions II}, 
\rm Ann. of Math. Studies, \bf{137}\rm, Princeton Univ. Press, Princeton, NJ (1995), 135-187.

\bibitem{FS4} J.E. Fornaess, N. Sibony, \it{Some open problems in higher dimensional complex analysis and complex dynamics}, \rm Publ. Mat., \bf{45} \rm (2001), no. 2, 529-547. 

\bibitem{HP} J.H. Hubbard, P. Papadopol, \it {Superattractive fixed points in ${\bf C}^n$}, 
\rm Indiana Univ. Math. J., \bf{43} \rm (1994), no. 1,  321-365. 



\bibitem{L1} F. Ledrappier, \it{Some properties of absolutely continuous invariant measure on an 
interval}, \rm Ergodic Theory Dynamical Systems, \bf{1} \rm (1981), no. 1, 77-93.

\bibitem{L2} F. Ledrappier, \it{Quelques propriétés ergodiques des applications rationnelles}, 
\rm C.R. Acad. Sci. Paris Sér. I Math., \bf{299} \rm (1984), no. 1, 37-40.

\bibitem{M} V. Mayer, \it{Comparing measures and invariant line fields}, \rm Ergodic Theory Dynamical 
Systems, \bf{22} \rm (2002), no. 2, 555-570. 


\bibitem{S} N. Sibony, \it {Dynamique des applications rationnelles de ${\bf P}^k$}, 
\rm in Dynamique et Géométrie Complexes, Panoramas et Synthèses No 8, SMF et EDP Sciences, 1999. 


\bibitem{Z} A. Zdunik, \it{Parabolic orbifolds and the dimension of the maximal measure for rational 
maps}, \rm Invent. Math., \bf{99} \rm (1990), no. 3, 627-649. 







\end{thebibliography}
\end{document}